\documentclass[a4paper,11pt]{article}
\RequirePackage{amsthm}
\RequirePackage{amssymb}
\usepackage[all,ps]{xy}
\usepackage{latexsym}
\newtheorem{thm}{Theorem}[section]
\newtheorem{prop}[thm]{Proposition}
\newtheorem{lem}[thm]{Lemma}
\newtheorem{cor}[thm]{Corollary}

\makeatletter
 \def\ps@notecras{%
    \def\@oddfoot{\small\@date\hfil\slshape\@pauthor\hfil\upshape\thepage}
    \def\@oddfoot{\basdepage}%
   \def\@oddhead{\tetedepage}
  \let\@mkboth\@gobbletwo
 \let\sectionmark\@gobble
\let\subsectionmark\@gobble}
\makeatother
\title{\bf Subvarieties of general type on a general projective hypersurface}
\author{Gianluca Pacienza}
\date{}
\begin{document}
\maketitle
%
%
%
\begin{quote}
{\small{\bf Abstract.} {We study subvarieties of a 
general projective degree $d$ hypersurface $X_d\subset \mathbf P^n$. 
Our main theorem, which improves previous results of L. Ein
and C. Voisin, 
implies in particular the following sharp corollary: any subvariety of
a general hypersurface $X_{d}\subset {\mathbf P}^n$, for
$n\geq 6$ and $d\geq 2n-2$, is of general type.}}
\end{quote}

%
%
\section {Introduction}
%
Let $X_d\subset \mathbf P^n$ be a general (in the countable Zariski
topology) complex projective hypersurface of degree $d$. 
The study of the geometry of $k$-dimensional subvarieties of $X_d$ 
in terms of $k,n$ and $d$ has received much attention in the last
15 years (see [C1], [E1] and [E2], [X1] and [X2], 
[V1] and [V2], [CL], [CLR],
[P], [C2], [CR]). In particular 
this study is related to the hyperbolicity of the hypersurface
$X_d\subset \mathbf P^n$. Recall that a compact complex manifold $M$ is 
said to be {\it hyperbolic}  (in the sense of Brody or Kobayashi) if 
there are no nonconstant entire holomorphic maps $f:{\mathbb C}\to M$.
S. Lang conjectured (cf. [L], Conjecture 5.6) that, in the case of
a projective variety $V$, the notion of hyperbolicity has an {\it algebraic}
characterization, namely $V$ is hyperbolic if and only if
any subvariety $Y$ of $V$ is of general type
(that is, if $Y$ is smooth, some multiple of the canonical bundle of $Y$
gives a projective embedding of a non-empty Zariski open subset of $Y$. If 
$Y$ is singular, then it is said to be of general type if 
some desingularization of $Y$ has this property). Notice that
if any subvariety $Y$ of $V$ is of general type, then in particular
$V$ does not contain rational curves or abelian subvarieties - a
condition which is of course implied by the hyperbolicity.
In this paper we focus our attention on the case of a 
general projective hypersurface $X_d\subset \mathbf P^n$, and 
give a sharp bound on its degree $d$ in order to satisfy the algebraic
property above that, conjecturally, should be equivalent to   
the hyperbolicity of $X_d$.  

This problem has been studied by L. Ein in [E1] and [E2], where, generalizing
a previous result by H. Clemens [C1], he proves
in particular that whenever $d\geq 2n-k$, for $n\geq 3$, then any 
$k$-dimensional subvariety $Y$ of the general $X_d\subset \mathbf P^n$ 
has nonzero geometric genus, and if the 
inequality is strict then $Y$ is of general type. Ein's result, which concerns 
more generally subvarieties of general complete intersections in an arbitrary 
smooth projective variety, has been improved by one, in the case 
of projective hypersurfaces, by C. Voisin ([V1], [V2]), who proves 
the following 

\bigskip

\noindent
{\bf Theorem (Voisin).}
 {\it Let $X_d\subset \mathbf P^n$ be a general hypersurface of degree 
 $d\geq 2n-k-1$, where $k$ is an integer such that $1\leq k\leq n-3$.
 Then any $k$-dimensional subvariety $Y$ of $X$
 has nonzero geometric genus, and if the 
 inequality is strict then $Y$ is of general type.}
\bigskip

Our main result is the following
\begin{thm}\label{mainthm}
 Let $X_d\subset \mathbf P^n$ be a general 
 complex projective hypersurface of degree $d$. Let $Y\subset X_d$ be a 
 $k$-dimensional subvariety with a desingularization 
 $\tilde Y{\buildrel j\over \longrightarrow} Y$ such that 
 $h^0 (\tilde Y,K_{\tilde Y}\otimes j^* {\cal O}_{\mathbf P^n}(-1))=0$.
 If the inequalities
 \begin{equation}\label{numcond1}
 d-1\geq \max \left\{ \frac{7n-3k-3}{4},\frac{3n-k+1}{2}\right\}
 \end{equation}
and 
 \begin{equation}\label{numcond2}
  {{d(d-3)}\over {2}}\geq 2n-k-3
 \end{equation}
are satisfied, then $Y$ is contained in the locus covered by the lines 
 of $X$.
\end{thm}

Since the dimension of the locus spanned by the lines on a general
hypersurface $X_d$ is equal
to $(2n-2)-d$, it follows from Theorem \ref{mainthm} that for 
any subvariety $Y$ of $X_d\subset \mathbf P^n$, whenever 
$d\geq 2n-2$ and $n\geq 6$, the canonical bundle of 
a desingularization $j:\tilde Y\to Y$ is the sum of the effective divisor 
$K_{\tilde Y}\otimes j^* {\cal O}_{\mathbf P^n}(-1)$ and of 
$j^* {\cal O}_{\mathbf P^n}(1)$, which is very ample on an open subset,
so we obtain
\begin{cor}\label{cor}
Any subvariety of a general hypersurface 
$X_{d}\subset {\mathbf P}^n$, with $d\geq 2n-2$ and
$n\geq 6$, is of general type.
\end{cor}

Corollary \ref{cor} is sharp, since general hypersurfaces of degree $2n-3$ 
contain a finite number of lines (and, by [P], lines are the only rational 
curves allowed on the general $X_{2n-3}\subset \mathbf P^n$, for $n\geq 6$).

The weird looking numerical hypotheses (1) and (2) 
of Theorem \ref{mainthm} are needed in order 
to control the positivity of the twisted exterior powers of the bundle
$M^d_{{\mathbf P}^n}$ (resp. $M^d_{G}$) over ${\mathbf P}^n$
(resp. over the Grassmannian of lines in ${\mathbf P}^n$), that are 
defined in \S 2. This control will appear to be a crucial point along
the proof (cf. Lemmas \ref{Clemlemma}, \ref{pos2} and \ref{pos3}).
The proof of Theorem \ref{mainthm} makes use of the powerful variational 
approach introduced by C. Voisin in [V1] and [V2], and adopted by the 
author in [P] to study the geometry of subvarieties having 
geometric genus zero on a general hypersurface. These methods have been
strengthten more recently by H. Clemens and Z. Ran (see [C2] and [CR]) 
to study in greater generality subvarieties
$Y$ of $X$ 
with desingularizations $j: \tilde Y\to Y$ verifying 
$h^0 (\tilde Y, K_{\tilde Y}\otimes j^* {\cal O}_{\mathbf P^n}(a))=0$,
for some integer $a\geq 0$. 

The proof is naturally divided in two parts.
First, following an idea that goes back to
Voisin  [V2], under a technical numerical hypothesis, 
it is possible to see that through each point of $Y$ there is a line which
intersects set-theoretically $X$ in at most two points.
Precisely we prove
\begin{prop}\label{main}
 Let $X_d\subset \mathbf P^n$ be a general hypersurface, whose 
 degree $d$ verifies the numerical condition (\ref{numcond1}) of 
 Theorem \ref{mainthm}, 
 and $Y\subset X_d$ a subvariety 
 of dimension $k$ such that 
 $h^0 (\tilde Y,K_{\tilde Y}\otimes j^* {\cal O}_{\mathbf P^n}(-1))=0$, 
 where $j:\tilde Y\to Y$ is
 a desingularization.
 Then, for some $r\geq 1$, $Y$ is contained in the 
 sublocus $\Delta_{(r,d-r),X}$ of $X_d$ defined as
 $$  
  \Delta_{(r,d-r),X}:=\lbrace x\in X_d : \exists\ a\ line\ \ell\ 
  s.t.\ \ell\cap X_d= r{\cdot} x+(d-r){\cdot} x',\ x'\in X_d\rbrace.
 $$
\end{prop}
(This result, under different numerical hypothesis, 
can also be found in [CR]).

The second part of the proof of Theorem \ref{mainthm}
deals with the study of the 
locus ${\Delta}_{(r,d-r),X}$. For this, 
I use two explicit desingularizations of
${\Delta}_{(r,d-r),X}$, 
whose canonical bundles are easily computable. Then, in both cases,
the key point is the construction of a globally generated 
subbundle contained in the exterior powers of the  
(twisted) tangent bundle of the family of the desingularizations. 
This fact will allow us to obtain
the following proposition which concludes the proof of Theorem 
\ref{mainthm}:

\begin{prop}\label{mainprop}
 Let $X_d\subset \mathbf P^n$ be a general hypersurface of degree $d$ verifying 
 the numerical condition (\ref{numcond2}) of 
 Theorem \ref{mainthm}.
 Let $Y\subset {\Delta}_{(r,d-r),X}$ 
 be a subvariety of dimension $k$, and $j:{\tilde Y}\to Y$ 
 a desingularization such that
 $h^0 (\tilde Y,K_{\tilde Y}\otimes j^* {\cal O}_{\mathbf P^n}(-1))=0$ .
 Then $Y$ is contained in the locus of lines of $X_d$.
\end{prop}

\noindent
{\small {\it Acknowledgements.} I would like to thank M. Zaidenberg 
and L. Manivel, who invited me, respectively on February 2001 and March 2002, 
at the Institut Fourier in Grenoble, where the preliminary version of
this work was started and brought to end.
My sincerest thanks to H. Clemens, who sent me the 
manuscript [C2] 
and kindly encouraged me when I communicated  
to him the idea of adapting some arguments contained in 
[C2] to strengthten my work.
Finally, it is a pleasure to acknowledge the great influence 
of C. Voisin's articles [V1] and [V2] on my work.}

%
\section {Preliminaries}
%
We will follow the notation already used in [P], which we recall below.

\noindent
{\bf Notation.}

\noindent $S^d:=H^0 (\mathbf P^n, {\cal  O}_{\mathbf P^n}(d))$;
\\ 
$S^d_x:=H^0 (\mathbf P^n, {\cal I}_x\otimes{\cal  O}_{\mathbf P^n}(d))$; 
\\ 
$N:=h^0(\mathbf P^n, {\cal  O}_{\mathbf P^n}(d))= dim\ S^d$;
\\ 
${\cal X}\subset {\mathbf P^n}\times S^d$ 
will denote the universal hypersurface of degree $d$;
\\
$X_F\subset {\mathbf P^n}$ the fiber of the family ${\cal X}$
over $F\in S^d$, i.e. the hypersurface 
defined by $F$. 
\smallskip

Let $U\rightarrow S^d$ be an \'etale map and
${\cal Y}\subset {\cal X}_U$ a reduced and 
irreducible subscheme of relative dimension $k$
(in the following, by abuse of notation, we will often omit
the \'etale base change).
Let $\tilde {\cal Y}\rightarrow {\cal Y}$ be a desingularization and 
$\tilde {\cal Y}{\buildrel j\over\longrightarrow}{\cal X}_U$ 
the natural induced map.
We may obviously assume $\cal Y$ invariant under some 
lift of the action of $GL(n+1)$ (recall that $g\in GL(n+1)$
acts on the product ${\mathbf P^n}\times S^d$ as follows:
$
g(x,F)= (g(x), (g^{-1})^*F)).
$
Let $\pi:{\cal X}\rightarrow {\mathbf P^n}$ be the projection on the
first component and  $T^{vert}_{\cal X}$ (resp. $T^{vert}_{\cal Y}$)
the vertical part of $T_{\cal X}$ (resp. $T_{\cal Y}$) 
w.r.t. $\pi$, i.e. $T^{vert}_{\cal X}$ (resp. $T^{vert}_{\cal Y}$) is the sheaf
defined by
$$
 0\rightarrow T^{vert}_{\cal X}\rightarrow T{\cal X}{\buildrel {\pi_*}\over
 \longrightarrow} T{\mathbf P^n}\rightarrow 0
$$
$$
 ({\it {resp.}}\ \ \ Ý 0\rightarrow T^{vert}_{\cal Y}\rightarrow 
 T{\cal Y}\ {\buildrel {\pi_*}\over \longrightarrow} T{\mathbf P^n}).
$$
\bigskip

\noindent
{\bf Remark 2.0.1.}
{\it Let $\cal Y$ be a subscheme of 
${\cal X}\subset {\mathbf P^n}\times S^d$ of
relative dimension $k$ and invariant under the action of $GL(n+1)$. Then:}
\begin {description}
\item[\it (i)]
{\it the map $T{\cal Y}\ {\buildrel {\pi_*}\over 
\longrightarrow} T{\mathbf P^n}$ 
is surjective and hence
$$
 codim_{\ T^{vert}_{{\cal X},(y,F)}}\ T^{vert}_{{\cal Y},(y,F)}
 =codim_{\cal X}{\cal Y}=n-k-1.
$$
}
\noindent

\item[\it (ii)]
{\it
$T^{vert}_{{\cal Y},(y,F)}$ contains the vertical part of the
tangent space to the orbit 
of the point $(y,F)$ under the action of $GL(n+1)$,
i.e. 
$$
 T^{vert}_{{\cal Y},(y,F)}\supset\ <S^1_y \cdot J_F^{d-1}, F>,
$$
where $J_F^{d-1}$ is the Jacobian ideal of $F$.}

\end {description}
\bigskip

Let $d$ be a positive integer.
Consider the bundle $M^d_{{\mathbf P}^n}$ defined by the exact sequence
\begin{eqnarray}\label{M_d}
 0\rightarrow M_{{\mathbf P}^n}^d
 \rightarrow S^d \otimes {\cal O}_{\mathbf P^n}
 {\buildrel {ev}\over
 \longrightarrow} {\cal O}_{\mathbf P^n}(d)\rightarrow 0,
\end{eqnarray}
whose fiber at a point $x$ identifies by definition to $S_x^d$. 
From the inclusion ${\cal X}\hookrightarrow {\mathbf P^n}\times S^d$ 
we get the exact sequence
$$
 0\rightarrow T{\cal X}_{|X_F}\rightarrow 
 T{\mathbf P^n}_{|X_F}\oplus (S^d\otimes {\cal O}_{X_F})
 \rightarrow {\cal O}_{X_F}(d)\rightarrow 0,
$$
which combined with (\ref{M_d}) gives us
\begin{eqnarray}\label{vert}
 0\rightarrow {M_{{\mathbf P}^n}^d}_{|X_F}\rightarrow T{\cal X}_{|X_F}
 \rightarrow T{\mathbf P^n}_{|X_F}\rightarrow 0.
\end{eqnarray}
In other words 
${M_{{\mathbf P}^n}^d}_{|X_F}$ identifies to the 
vertical part of $T_{\cal X}\otimes {\cal O}_{X_F}$ with respect to the 
projection to ${\mathbf P^n}$.
 
Let $G:=Grass(1,n)$ be the Grassmannian of lines in ${\mathbf P}^n$,
${\cal O}_G (1)$ the line bundle on $G$ giving its Pl\"ucker polarization,
and ${\cal E}_d$ be the $d^{th}$-symmetric power of the dual of the 
tautological subbundle on $G$. 
Recall that the fibre of ${\cal E}_d$ at a point $[{\ell}]$ 
is, by definition, given by $H^0 ({\ell},{\cal O}_{\ell} (d))$.

Let $M_G^d$ be the vector bundle on $G$ 
defined as the kernel of the evaluation map:
$$
 0\rightarrow M_G^d\rightarrow S^d \otimes {\cal O}_G 
 \rightarrow {\cal E}_d \rightarrow 0.
$$
Notice that the fiber of $M_G^d$ at a point 
$[\ell]$ is equal to $H^0 ({\cal I}_{\ell}(d))$. 

The bundles $M_{{\mathbf P}^n}^d$ and $M_G^d$ verify the following 
positivity properties, that will be often used in what follows:

\begin{lem}\label{pos}
 (i) 
 $M_{{\mathbf P}^n}^d \otimes {\cal O}_{{\mathbf P}^n}(1)$ is 
 generated by its global sections;\\
 (ii) $M_G^d\otimes {\cal O}_G (1)$ is 
 generated be its global sections.
\end{lem}
\begin{proof} See, for instance, [P].

\end{proof}

%
\section {A first reduction}
%

Let
$Y_F\subset X_F$ be a general ($k$-dimensional) fiber of the subfamily 
${\cal Y}\subset {\cal X}_U$, and 
$\tilde Y_F{\buildrel j\over \longrightarrow} Y_F$ its desingularization. 
By abuse of notation, we will often write 
$K_{\tilde Y_F}(-1)$ 
instead of
$K_{\tilde Y_F}\otimes j^* {\cal O}_{\mathbf P^n}(-1))$.

Recall now the isomorphisims:
\smallskip

(i) ${\Omega ^{N+k}_{{\tilde{\cal Y}}}}_{|{\tilde Y_F}}
\cong K_{\tilde Y_F}$;
\smallskip

(ii) $(\wedge^{n-1-k} T{{\cal X}_U}_{|{X_F}})\otimes K_{X_F}
\cong {\Omega ^{N+k}_{{\cal X}_U}}_{|{X_F}}$.
\smallskip

\noindent 
Set $c:={n-1-k}= codim_{X_F}\ Y_F$.
Using (i) and (ii), from the natural morphism 
$\Omega ^1_{{{\cal X}_U}}\rightarrow \Omega ^1_{{\tilde{\cal Y}}}$, we
get a map
\begin{equation}
 (\wedge^{c} T{{\cal X}_U}_{|{X_F}})\otimes K_{X_F}
 \cong {\Omega ^{N+k}_{{\cal X}_U}}_{|{X_F}}
 \rightarrow {\Omega ^{N+k}_{{\tilde{\cal Y}}}}_{|{\tilde Y_F}}
 \cong K_{\tilde Y_F},
\end{equation}

and hence, after tensoring by ${\cal O}_{X_F}(-1)$,
\begin{equation}\label{restriction}
 (\wedge^{c} T{{\cal X}_U}_{|{X_F}})\otimes K_{X_F}(-1)
 \longrightarrow K_{\tilde Y_F} (-1).
\end{equation}

Now taking global sections in (\ref{restriction}) and using (\ref{vert})
we have the following commutative diagram
\begin{eqnarray}\label{composite}
 &\ \ \ 
 H^0 (\wedge^{c} T{{\cal X}_U}_{|{X_F}}\otimes K_{X_F}(-1))&
 \longrightarrow   \ 
 H^0 (K_{\tilde Y_F}(-1))\nonumber\\
 &\cup 
 &\ \ \nearrow\nonumber\\
 &H^0 (\wedge^{c} {M_{{\mathbf P}^n}^d}_{|X_F}\otimes K_{X_F}(-1))
 &
\end{eqnarray}
 
By hypothesis, we have that the 
composite map in (\ref{composite}) is identically zero.
This, by the $GL(n+1)$-invariance of ${\cal Y}$, implies  
that $T^{vert}_{{\cal Y},(y,F)}$ is then contained 
in the base locus of  
$H^0 (\wedge^{c} {M_{{\mathbf P}^n}^d}_{|{X_F}} \otimes K_{X_F}(-1))$,
considered as the
space of sections of a line bundle on the Grassmannian of
codimension $c$ subspaces of $T^{vert}{\cal X}_{|{X_F}}$.

The generalization presented in [C2] of the variational approach
introduced by C. Voisin in [V1] and [V2], and applied by the 
author in [P], starts with a sharp 
algebraic study of the base locus of the bundles 
$\wedge^{c} {M_{{\mathbf P}^n}^d} (b)$. Precisely, we will use the following

\begin{lem}\label{Clemlemma}{\bf ([C2])}.
 Suppose $d$ verifies the numerical condition (\ref{numcond1}) 
 of Theorem \ref{mainthm}. 
 Let $A$ be a codimension $c$ subspace of $S^d_x=(T_{{\cal X},(x,F)})^{vert}$ 
 which is in the base locus of 
 $H^0 (\bigwedge^c {M_{{\mathbf P}^n}^d}(d-n-2))$. 
 Then there exists a line ${\ell_A}$ passing through $x$ 
 such that 
 \begin{eqnarray}\label{rankcond}
  rk\ {A \over A\cap H^0 ({\mathcal I}_{\ell_A} (d))}\leq n+1.
 \end{eqnarray}
\end{lem}
\begin{proof}
In [C2] there is the study of the base locus of 
$H^0 (\bigwedge^c {M_{{\mathbf P}^n}^d}(d-n-1+a))$, for $a\geq 0$ 
(this is the point where Lemma \ref{pos}, (i) is used). 
Here we simply remark that 
the arguments presented in [C2] also apply to the $a=-1$ case.
For the reader's convenience we outline the idea of the proof (for the
details see [C2]).
The main point in [V2] and [P] was to produce, 
by Koszul complexes' techniques,
{\it explicit} global sections of the bundle
$\wedge^{c} {M_{{\mathbf P}^n}^d} (b)$, for the integers
$b,c,d$ considered there. This was used to deduce that,
for a generic polynomial  $P\in S^{d-1}$,
the multiplication map 
\begin{eqnarray}
  m_{P,A} : &S^1_x& 
  \rightarrow S^d_x/A \nonumber
  \\&L&\mapsto L\cdot P
  \ mod\ A\nonumber
\end{eqnarray}
has rank one.
H. Clemens considers more generally in [C2] the smallest integer $s\geq 0$
such that rank of the 
multiplication map
\begin{eqnarray}
  m_{P,A,s} : &S^1_x& 
  \rightarrow S^d_x/(A+ Q_1\cdot S^1_x+\ldots  + Q_s\cdot S^1_x) \nonumber
  \\&L&\mapsto L\cdot P
  \ mod\ (A+ Q_1\cdot S^1_x+\ldots  + Q_s\cdot S^1_x)\nonumber
  \\&&P, Q_1,\ldots ,Q_s \ are\ generic\ polynomials\ in\ S^{d-1}\nonumber
\end{eqnarray}
is one.
Then, as in [V2] and [P], an infinitesimal argument applies.
Namely, recall that if $V$ and $W$ are vector spaces, and 
 $Z_k:=\lbrace \phi\in Hom(V,W): rank\ \phi\leq k\rbrace$,
 then 
 \begin{equation}\label{tg}
  T_{{Z_k},\phi}=\lbrace \psi\in Hom(V,W): 
  \psi (ker\phi)\subset Im\phi\rbrace.
 \end{equation}
 Applying this to the map $m_{P,A,s}$, 
 we obtain that, for any $R\in S^{d-1}$,
 $$
 R\cdot Ker\ {m_{P,A,s}}
 \ mod\ (A+ Q_1\cdot S^1_x+\ldots + Q_s\cdot S^1_x)
 \subset Im\ m_{P,A,s} 
 $$ 
 i.e. 
 $$
  H^0 ({\mathcal I}_{\ell_A} (d))
  \subset A+ Q_1\cdot S^1_x+\ldots  + Q_s\cdot S^1_x + P\cdot S^1_x ,
 $$ where
 ${\ell_A}$ is the line determined by $Ker\ {m_{P,A,s}}$. 
To complete the proof it remains to verify that the line is 
independent of the choice 
of the polynomials, and that, under the hypothesis 
(\ref{numcond1}), the integer $s$ is such that $A$ and 
$H^0 ({\mathcal I}_{\ell_A} (d))$ satisfy (\ref{rankcond}).
\end{proof}
Since the map
(\ref{composite}) vanishes, Lemma \ref{Clemlemma} applies to 
the tangent space $T^{vert}_{{\mathcal Y},(y,F)}$.
Then, at a generic point $(y,F)\in {\mathcal Y}$,
the tangent space $T^{vert}_{{\mathcal Y},(y,F)}$ contains a subspace
$T\subset H^0 ({\mathcal I}_{\ell_{(y,F)}} (d))$, where $\ell_{(y,F)}$ 
is a line through 
$y$ and $T$ satisfies 
\begin{eqnarray}\label{rkcond2}
rk\ {T^{vert}_{{\cal Y},(y,F)} \over T}\leq n+1.
\end{eqnarray}
We now verify an easy fact:
\begin{lem}
The tangent space $T^{vert}_{{\mathcal Y},(y,F)}$
cannot contain two subspaces $T$ and $T'$ of  ideals of different
lines $\ell\not = \ell'$ and verifying (\ref{rkcond2}).
\end{lem}
\begin{proof} 
Indeed, if this were the case, by the surjectivity of
$$
 H^0 ({\mathcal I}_{\ell} (d))\oplus H^0 ({\mathcal I}_{\ell'} (d))
 \twoheadrightarrow S^d_y,
$$
and the numerical condition (\ref{rkcond2}),
then $T^{vert}_{{\mathcal Y},(y,F)}$ would contain a 
subspace of $S^d_y$ of codimension at most $2(n+c+1-d)$. Now, by remark 
{2.0.1}, (i), we have 
$$
 c=codim_{\ T^{vert}_{{\cal X},(y,F)}}\ T^{vert}_{{\cal Y},(y,F)}
 =codim_{\cal X}{\cal Y}\leq 2(n+c+1-d)
$$
which is equivalent to 
$$
 d\leq {3n-k\over 2},
$$
and the last inequality is impossible because of (\ref{numcond1}).
\end{proof}
Then, we can consider 
the distribution 
${\mathcal T}\subset T^{vert}_{{\mathcal Y}}$, 
pointwise given by the  $T$'s.
This distribution
turns out to have the following properties
\begin{prop}\label{int}
The distribution ${\mathcal T}\subset T^{vert}_{{\mathcal Y}}$ is
integrable and the natural 
map $\phi:{\mathcal Y}\to G(1,n)$, associating to $(y,F)$ 
the line determined by $T\subset H^0 ({\mathcal I}_{\ell_{(y,F)}} (d))$,
is constant along the leaves
of the corresponding foliation. 
\end{prop}
\begin{proof}
 The proof goes along the lines of [V2], lemma 3 and 4, and of [P], lemma 3.3.
 For the detailed proof in the general case, see [C2]. 
 Again, for the reader's convenience, we sketch it below.
 Consider the bracket map
  $$
  {\Psi}:\bigwedge^2 {\cal T}\rightarrow
  T^{vert}_{\cal Y}/{\cal T}\subset 
  {T^{vert}_{\cal X}}_{|{\cal Y}}/{\cal T},
 $$
 which is given at the point $(y,F)$ by 
 $$
  \psi : \wedge^2 T_{\ell _{(y,F)}}
  \rightarrow T^{vert}_{{\cal Y},(y,F)}/T_{\ell _{(y,F)}}
  \hookrightarrow H^0 ({\cal O}_{{\ell}_{(y,F)}} (d) (-y)).
 $$
 Now, choose coordinates on $\mathbf P^n$ such that 
 $\ell:= \ell _{(y,F)}
 =\lbrace X_2=\ldots =X_n =0\rbrace$ and $y=[1,0,\ldots ,0]$.
 Note that, since $y\in \ell$,  $\phi_* (T^{vert}_{{\cal Y},(y,F)})$
 is contained in $H^0 (N_{\ell/{\mathbf P^n}}(-y))$ 
 ($\phi_*$ is the differential of $\phi$ at ${(y,F)}$).
 One verifies that 
 \begin{eqnarray}\label{star}
  \psi (A\wedge B)= A\cdot \phi_* (B)-B\cdot \phi_* (A),
  \ A,B\in T_{\ell _{(y,F)}},
 \end{eqnarray}
 where the 
 bilinear map $(a,b)\mapsto a\cdot b$ is explicitely given by 
 $$ 
  P\cdot (X_1 \sum_{i=2}^n b_i {\partial\over{\partial X_i}}) = 
  \sum_{i=2}^n b_i X_1({\partial P\over{\partial X_i}})_{|\ell} 
  \in H^0 ({\cal O}_{{\ell}} (d) (-y)).
 $$
 A key linear algebra lemma allows to prove that $\phi_*$ is zero,
 so the proposition follows from (\ref{star}) and from the Frobenius theorem.
\end{proof}

Using this we will prove, via the $GL(n+1)$-invariance of ${\mathcal Y}$,
that $F_{|{\ell_{(y,F)}}}$ has, set-theoretically, at most two zeroes.

\smallskip
\noindent {\it Proof of Proposition \ref{main}.} 
Let $(y,F)$ be a general point
of ${\mathcal Y}$. Let $\ell$ be the line through $y$ and 
$T$ the subspace of $H^0 ({\mathcal I}_{\ell} (d))$ contained in
$T^{vert}_{{\cal Y},(y,F)}$ and verifying (\ref{rkcond2}). 
Consider the
following diagram
\begin{equation}\label{codim}
\xymatrix{ 
 &0\ar[d]
&0\ar[d]
\\
0\ar[r]
&T\ar[d]\ar[r]
&H^0 ({\mathcal I}_{\ell} (d))\ar[d]
\\
0\ar[r]
&T^{vert}_{{\cal Y},(y,F)}\ar[d]\ar[r]
&T^{vert}_{{\cal X},(y,F)}=S^d_y\ar[d]
\\
 &H^0 ({\cal O}_{{\ell}} (d)(-y))\ar@{=}[r]
&H^0 ({\cal O}_{{\ell}} (d)(-y))\ar[d]
\\
 &
 &0
}\end{equation}
By (\ref{rkcond2}) we have
\begin{eqnarray}\label{dimim}
 dim\ Im\ (T^{vert}_{{\cal Y},(y,F)}
 \rightarrow H^0 ({\cal O}_{{\ell}} (d)(-y)))\leq n+1
\end{eqnarray}
 
On the other hand, by remark 2.0.1, (ii), the vertical tangent space 
$T^{vert}_{{\cal Y},(y,F)}$ contains
$S^1_y \cdot J_F^{d-1}$ and $F$ itself. 
Take coordinates $X_0,\ldots ,X_n$ on
$\mathbf P^n$ such that $y=[1,0,\ldots ,0]$, and
${\ell}:={\ell}_{(y,F)}=\lbrace X_2=\ldots =X_n=0\rbrace$.
Since $\phi$ is constant along the leaves of the foliation, 
we can generically choose a polynomial $G$ in the leaf integrating $T$ 
so that the $(n-1)$-elements
$X_1{\partial G\over{\partial X_i}},\ i\geq 2,$ are independent modulo
the subpace
$$
 K:= <G_{|{\ell}},
 X_1({\partial G\over{\partial X_0}})_{|{\ell}},
 X_1({\partial G\over{\partial X_1}})_{|{\ell}}>\ 
 \subset H^0 ({\ell}, {\cal O}_{{\ell}}(d)\otimes
 {\cal I}_y),
$$
which is uniquely determined by $F_{|{\ell}}$ and hence is
constant along the leaf integrating $T$.
By (\ref{dimim}), this implies $dim K\leq 2$, that is
$$
 F_{|{\ell}}=\alpha X_1^r L^{d-r},
$$
for some $r\geq 1$, and some linear form $L$ on $\ell$.
\hfill $\Box$

We are then led to study  the locus 
$\Delta_{(r,d-r),F}$. 
This will be done in the last section.

%
\section {The bicontact locus $\Delta_{(r,d-r),F}$}
%
Let $X_F\subset \mathbf P^n$ be a general 
hypersurface of degree $d$ verifying (\ref{numcond1}), and $Y_F\subset X_F$ a
$k$-dimensional subvariety  whose desingularization 
$\tilde Y$ is such that $h^0 (\tilde Y,K_{\tilde Y}(-1))=0$. 
Then, by Proposition \ref{main}, we know that
$Y_F$ is contained in $\Delta_{(r,d-r),X_F}\subset X_F$, 
the $(2n-d)$-dimensional subvariety
of points $x$ of $X_F$ through which there is an osculating line $\ell$
intersecting $X_F$ at most at another point, $\it i.e.$ 
$\ell\cap X_F= r{\cdot} x+(d-r){\cdot} x',\ x'\in X_F$.
In what follows we will write $\Delta_{(r,d-r),F}$
instead of $\Delta_{(r,d-r),X_F}$.
To prove our theorem, we study
two  explicit desingularizations
of ${\Delta}_{(r,d-r),F}$, 
which have been used in [V2],
that are both given in terms of the zero locus 
of a section of a vector bundle.
Thus we compute, by adjunction,
the canonical bundle of such a desingularization.
Then, again, we adopt a variational approach and 
construct, in both cases, a  
subbundle contained in the exterior powers of the
(twisted) tangent bundle to the family of the desingularizations. 
A positivity result, namely the global generation of this subbundle, 
allows us to conclude the proof.

\bigskip

\noindent
{\bf Case 1: $r\geq 2$ and $d-r\geq 2$.}

\noindent
Let $G:= Gr(1,n)$ be the Grassmannian of lines in $\mathbf P^n$. 
Let ${\cal O}_G (1)$ be the line bundle on $G$ which gives the Pl\"ucker 
embedding.
Let $Z$ be the blow-up along the diagonal $\Delta$ of the product 
${\mathbf P}^n\times {\mathbf P}^n$ with projections:
\begin{equation}
\xymatrix{
{Z:= Bl_{\Delta}{\mathbf P}^n\times {\mathbf P}^n}\ar[r]^{\ \ \ \ \ b}
& {\mathbf P}^n \times{\mathbf P}^n\ar[r]^{p_2}\ar[d]^{p_1}
& \mathbf P^n.\\
&\mathbf P^n
& 
& \\
}\end{equation}

Consider the map 
\begin{equation}
f\ :\ Z \to Gr(1,n),
\end{equation}
$$
 z\mapsto \ell_{z}
$$
where $\ell_z$ is the line determined by $z$.
Let $\tilde{p_i}:= p_i\circ b$, for $i=1,2$, and consider the line bundles
on $Z$ defined as follows:
$H_i:=\tilde{p_i}^* {\cal O}_{\mathbf P^n}(1)$ and $L:=f^* {\cal O}_G (1)$.
The variety $Z$ comes together with a projective bundle:
${\mathbf P}{\buildrel \pi \over {\longrightarrow}} Z$, and we define  
${\cal E}_d:= \pi_* {\mathcal O}_{\mathbf P} (d)$. Notice that the fibre 
of ${\cal E}_d$ at $z$ is equal to $H^0 ({\ell_z},{\cal O}_{\ell_z} (d))$.
Consider the line bundle ${\mathcal L}_{r,d-r}\subset {\mathcal E}_d$,
whose fibre at $z\in Z$ is given by the one dimensional space of
polynomials $P\in H^0 ({\ell_z},{\cal O}_{\ell_z} (d))$ vanishing at $x$ to 
the order $r$ and at $y$ to the order $d-r$, where 
$(x,y)=b(z)\in {\mathbf P}^n \times{\mathbf P}^n$. 
Define ${\mathcal F}_{r,d-r}:={\mathcal E}_d/{\mathcal L}_{r,d-r}$. 
To any polynomial 
$F\in S^d$ we can associate a section $\sigma_F\in H^0 (Z,{\mathcal E}_d)$, 
whose value at a point 
$z$ is exactly the polynomial $F_{|\ell_z}\in {{\mathcal E}_d}_{|z}$, 
and we will denote 
by ${\bar{\sigma}_F}$ its image 
in $H^0 (Z,{\mathcal F}_{r,d-r})$. Then we define 
${\tilde\Delta}_{(r,d-r),F}:= V({\bar{\sigma}_F})$. 
By construction, we have 
$\tilde{p_1}({\tilde\Delta}_{(r,d-r),F})={\Delta}_{(r,d-r),F}$.
Since ${\mathcal F}_{r,d-r}$ is generated by its global sections, the variety
${\tilde\Delta}_{(r,d-r),F}$ is smooth and of the right dimension.
Moreover, since in the degree considered through a generic 
point of ${\Delta}_{(r,d-r),F}$ there is just one $r$-osculating line,
the map $\tilde{p_1}$ is a desingularization of ${\Delta}_{(r,d-r),F}$ .
We will now recall how to compute the canonical bundle of 
${\tilde\Delta}_{(r,d-r),F}$.
As remarked in [V2], the Picard group of $Z$ is generated by 
$H_1,\ H_2$ and $L$, the canonical class of $Z$ is 
$K_Z= -2H_1 - 2H_2 +(-n+1)L$, and the class of ${\mathcal L}_{(r,d-r),F}$
is given by $rH_1 + (d-r)H_2$.
Therefore, by adjunction, we have
$$
 K_{{\tilde\Delta}_{(r,d-r),F}}= K_Z + c_1 ({\mathcal F}_{r,d-r})=
 (r-2)H_1 + (d-r-2)H_2 + ( {{d(d-1)}\over {2}} - n+1 ) L.
$$

Consider now the bundles ${\mathcal N}^{r,d-r}_Z$ and ${\mathcal M}^{d}_Z$ 
on $Z$ respectively 
defined by the two following exact sequences :
\begin{eqnarray}\label{N_(r,d-r)}
 0\rightarrow {\mathcal N}^{r,d-r}_Z
 \rightarrow S^d \otimes {\mathcal O}_{Z}
 \rightarrow {\mathcal F}_{r,d-r}\rightarrow 0,
\end{eqnarray}

\begin{eqnarray}\label{M}
 0\rightarrow {\mathcal M}^{d}_Z
 \rightarrow S^d \otimes {\mathcal O}_{Z}
 \rightarrow {\mathcal E}_{d}\rightarrow 0.
\end{eqnarray}

By definition we have
\begin{eqnarray}\label{inclusion}
 0\rightarrow {\mathcal M}^{d}_Z
 \rightarrow {\mathcal N}^{r,d-r}_Z
 \rightarrow {\mathcal L}_{r,d-r}\rightarrow 0.
\end{eqnarray}

The needed positivity result is the following
\begin{lem}\label{pos2}
 If $r\geq 3$, $d-r\geq 2$ and 
 \begin{equation}\label{cond2}
  \frac{d\left( d-1\right) }{2} -n\geq c -1, 
 \end{equation}
 then  
 the bundle 
 ${\wedge^c{\mathcal M}^{d}_Z}_{|{\tilde\Delta}_{(r,d-r),F}}\otimes 
 K_{{\tilde\Delta}_{(r,d-r),F}} (-H_1)$ is generated by its global sections.
\end{lem}
\begin{proof}
 Remark that ${{\mathcal M}^{d}_Z}=f^*{{M}^{d}_G}$.
 Hence, by lemma \ref{pos}, (ii), the bundle
 \begin{eqnarray}
  &&\wedge^c{\mathcal M}^d_Z \otimes det {\cal F}_{r,d-r}\otimes K_{Z}(-H_1)=\nonumber\\
  \!\!\! &=&\!\!\!\! f^* (\wedge^c M^d_G) \otimes {\cal O}_{Z}((r-3)H_1+ 
  (d-r-2)H_2 + ({d(d-1)\over 2}-n+1)L)\nonumber\\
  &=& \!\!\!\! f^* (\wedge^c M^d_G (c)) 
 \otimes {\cal O}_{Z}((r-3)H_1+(d-r-2)H_2 + 
  ({d(d-1)\over 2}-n -c+1)L),\nonumber
 \end{eqnarray}
 is globally generated under our numerical hypothesis, 
 and the same holds for its
 restriction to ${\tilde\Delta}_{(r,d-r),F}$.
\end{proof}

Let now $\Delta_{r,d-r}\subset {\mathbf P}^n\times S^d$ be the family of the 
$\Delta_{(r,d-r),F}$'s,
and $\tilde\Delta_{r,d-r}\subset {Z}\times S^d$ the family of the 
desingularizations. Let ${\cal Y}\subset \tilde\Delta_{r,d-r}$ be a subscheme
of relative dimension $k$, invariant under the action of
$GL(n+1)$, and $\tilde {\cal Y}\rightarrow {\cal Y}$ a desingularization.
Assume $h^0 ({\tilde Y}_F,K_{{\tilde Y}_F} (-H_1))=0$ and set $c=n-1-k$.
Recall the isomorphisms 
\begin{eqnarray}
 \wedge^{c}{{T{\tilde \Delta_{r,d-r}}}}_{|{{\tilde\Delta}_{(r,d-r),F}}}
 \otimes K_{\tilde\Delta_{(r,d-r),F}}&\cong&
 {{\Omega}^{N+k}_{{\tilde \Delta_{r,d-r}}}}_{|\tilde \Delta_{(r,d-r),F}}\\ 
 {\Omega ^{N+k}_{{\tilde{\cal Y}}}}_{|{\tilde Y_F}}
 &\cong& K_{\tilde Y_F}
\end{eqnarray}
and consider the natural map 
\begin{eqnarray}\label{untwist}
 \wedge^{c} {{T{\tilde \Delta_{r,d-r}}}}_{|{{\tilde\Delta}_{(r,d-r),F}}}
 \!\!\!\otimes K_{\tilde\Delta_{(r,d-r),F}}\!\!\!\!\!\cong\!
 {{\Omega}^{N+k}_{{\tilde \Delta_{r,d}}}}_{|\tilde \Delta_{(r,d-r),F}}\!\!
 \to
 {\Omega ^{N+k}_{{\tilde{\cal Y}}}}_{|{\tilde Y_F}}\!\!\!\!\!\!
 &\cong&\!\!\!\! K_{\tilde Y_F}.
\end{eqnarray}
If we twist (\ref{untwist}) by $-H_1$, then, by assumption, 
the induced map in cohomology 
\begin{equation}\label{S-restriction}
 H^0 (\wedge^{c}{{T{\tilde \Delta_{r,d-r}}}}_{|{{\tilde\Delta}_{(r,d-r),F}}}
 \otimes K_{\tilde\Delta_{(r,d-r),F}} (-H_1))
 \rightarrow H^0 (K_{{\tilde Y}_F} (-H_1))
\end{equation}
is zero. Let $T^{vert}_{\tilde \Delta_{r,d-r}}$ be the sheaf defined by
$$
 0\rightarrow T^{vert}_{\tilde \Delta_{r,d-r}}
 \rightarrow T{\tilde \Delta_{r,d-r}}
 {\rightarrow} T{Z}
 \rightarrow 0.
$$
Its restriction to ${{\tilde\Delta}_{(r,d-r),F}}$ coincides with 
${{\mathcal N}^{r,d-r}_Z}_{|{\tilde \Delta_{(r,d-r),F}}}$.
Therefore, by (\ref{inclusion}) and lemma \ref{pos2}, we have 
constructed a sub-bundle  
$$
 \wedge^c{{\mathcal M}^{d}_Z}_{|{\tilde\Delta}_{(r,d-r),F}}\!\!\otimes 
 K_{{\tilde\Delta}_{(r,d-r),F}} (-H_1)\!\!\hookrightarrow
 \wedge^c{{T^{vert}{\tilde \Delta_{r,d-r}}}}_{|{{\tilde\Delta}_{(r,d-r),F}}}
 \!\!\otimes K_{\tilde\Delta_{(r,d-r),F}} (-H_1),
$$ 
which is generated by its global sections, under the numerical hypothesis of 
the lemma \ref{pos2}.

We conclude the case 1 with the following
\begin{prop}\label{case1}
 Let $F$ be a general polynomial of degree $d$ verifying (\ref{cond2}),
 with $c=n-1-k$.
 Suppose $r\geq 2$ and $d-r\geq 2$.
 Let $Y_F\subset {\tilde\Delta}_{(r,d-r),F}$ 
 be a subvariety of dimension $k$, and $j:{\tilde Y}_F\to Y_F$ 
 a desingularization such that 
 $h^0 ({\tilde Y}_F,K_{{\tilde Y}_F}(-j^*H_1))=0$.
 Then $Y_F$ is contained in the locus of lines of $X_F$.
\end{prop}
\begin{proof}
When no confusion is possible, we will omit in what follows the 
index $(r, d-r)$ and simply set $\Delta=\Delta_{r,d-r}$.  
Suppose first $r\geq 3$.
Let $W\subset {T_{\tilde \Delta}}_{,(z,F)}$ 
 be a codimension $c$ subspace contained 
 in the base locus of 
 $H^0 (\wedge^c{{T{\tilde \Delta}}}_{|{{\tilde\Delta}_{F}}}
 \otimes K_{\tilde\Delta_{F}}(-H_1))$, considered 
 as the space of sections of a line bundle
 on the Grassmannian of codimension $c$ subspaces of 
 ${{T{\tilde \Delta}}}_{|{{\tilde\Delta}_{F}}}$. 
 Then we must have
 \begin{eqnarray}\label{vert.bl}
  W^{vert}:=W\cap {\mathcal N}^{r,d-r}_{Z}|_z \subset {\mathcal M}^{d}_{Z}|_z.
 \end{eqnarray}
 Indeed, if this were not the case, we would have 
 $codim_{{\mathcal M}^d_{Z}|_z} \bar W =c$, where 
 $\bar W:= W\cap {\mathcal M}^d_{Z}|_z$.
 Then consider the following commutative diagram:
\begin{equation}\label{finaldiagram}
  \xymatrix{
   H^0 (\wedge^c{{\mathcal M}^{d}_Z}_{|{\tilde\Delta}_{F}}\otimes 
   K_{{\tilde\Delta}_{F}} (-H_1))
   \ar[d]^{ev\ \ }\ar@^{(->}[r]\ar@_{>}[r]
   &H^0 (\wedge^c{{T{\tilde \Delta}}}_{|{{\tilde\Delta}_{F}}}
   \otimes K_{\tilde\Delta_{F}}(-H_1))\ar[d]^{ev\ \ }\\
   (\wedge^c{{\mathcal M}^d_Z}_{|{\tilde\Delta}_{F}}\otimes 
   K_{{\tilde\Delta}_{F}} (-H_1))|_z
   \ar[dr]^{<\cdot,{\bar W}>}\ar@^{(->}[r]\ar@_{>}[r]
   &(\wedge^c{{T{\tilde \Delta}}}_{|{{\tilde\Delta}_{F}}}
   \otimes K_{\tilde\Delta_{F}}(-H_1))|_z
   \ar[d]^{<\cdot,W>}\\
   &{\mathbb C}\\
}\end{equation}
($ev$ is the evaluation of the sections at the point $z$, and
 $<\cdot,W>$ is the contraction defined by the subspace $W$).
 Since $W$ belongs to the base locus of 
 $H^0 (\wedge^c{{T{\tilde \Delta}}}_{|{{\tilde\Delta}_{F}}}
 \otimes K_{\tilde\Delta_{F}}(-H_1))$, then the composite map 
 $<\cdot,{W}>\circ\ ev$ is zero, and 
 so would be $<\cdot,{\bar W}>\circ\ ev$. But this is absurd,
 because, by Lemma \ref{pos2}, the bundle
 $\wedge^c{{\mathcal M}^{d}_Z}_{|{\tilde\Delta}_{F}}\otimes 
 K_{{\tilde\Delta}_{F}} (-H_1)$
 is generated by its global sections.
 
 Let then ${\cal Y}\subset {\tilde\Delta}$ be a subvariety, 
 which is stable under the action of $GL(n+1)$ and of relative 
 codimension $c$. Assume moreover 
 that the restriction map (\ref{S-restriction})
 is zero.
 By (\ref{vert.bl}), $T^{vert}_{{\cal Y},(z,F)}$ is contained in
 \begin{equation}\label{etoile}
 {{\mathcal M}^d_Z}_{|z}=\lbrace G\in S^d : G_{|\ell_z}=0\rbrace.
 \end{equation} 
 On the other hand, by Remark 2.0.1, (ii), 
 $T^{vert}_{{\cal Y},(z,F)}$ contains  
 $F$ itself.
 So by (\ref{etoile}) we have that $F_{|\ell_z}=0$ 
 for every point $z\in Y_F$, 
 i.e. $Y_F$ is contained in the subvariety covered
 by the lines contained in $X_F$. 

 If $r=2$, 
 we can consider the natural isomorphism
 \begin{equation}\label{keyiso}
  {\tilde\Delta}_{(r,d-r),F}\ \ 
  {\tilde{\longrightarrow}}\ {\tilde\Delta}_{(d-r,r),F}
 \end{equation}
 sending a point $z\in {\tilde\Delta}_{(r,d-r),F}$, with $b(z)=(x,x')$, 
 to a  point 
 $w\in {\tilde\Delta}_{(d-r,r),F}$ with $b(w)=(x',x)$, where $x$ and $x'$
 are points on $X_F$ linked by the condition
 $$
  \exists\ a\ line\ \ell\ s.t.\ 
 \ell\cap X_F= r{\cdot} x+(d-r){\cdot} x'.
 $$
 Since $r=2$ implies $d-r\geq 3$, from what we have done before 
 it follows that the conclusion is true for
 ${\tilde\Delta}_{(d-r,r),F}$, and so, by (\ref{keyiso}), 
 the same holds for ${\tilde\Delta}_{(r,d-r),F}$.
 \end{proof}

\bigskip

\noindent
{\bf Case 2: $r=1$ or $d-r=1$.}

\noindent
Suppose for instance $d-r=1$. 
Let $\Gamma\subset \mathbf P^n\times  Gr(1,n)$ be the 
incidence variety, and $p$ and $q$ the projections on the two factors.
Let $\pi:{\mathbf P}\to \Gamma$ be the pull-back of the universal 
$\mathbf P^1$-bundle over $Gr(1,n)$ and $\tau$ the natural section of $\pi$.
Consider the bundle ${\mathcal E}_d:= \pi_* {\mathcal O}_{\mathbf P} (d)$  
over $\Gamma$, and its rank $2$ subbbundle 
${\mathcal K}\subset {\mathcal E}_d$
such that its fiber at a point $(x,\ell)$ is given by the
polynomials $P\in H^0 ({\mathcal O}_{\ell}(d))$ vanishing to the order 
at least $(d-1)$ at $x$. 
Consider the line bundle ${\mathcal L}_1$ defined by 
$$
 0\to {\mathcal L}_1\to {\mathcal E}_1\to 
 \tau^* {\mathcal O}_{\mathbf P} (1)=:H\to 0.
$$
Then ${\mathcal K}\cong {\mathcal L}_1^{d-1}\otimes {\mathcal E}_1$. 
Let ${\mathcal F}_d$ be the quotient ${\mathcal E}_d/{\mathcal K}$.
As in Case 1, to any $F\in S^d$ we can associate a global 
section $\sigma_F$ of ${\mathcal F}_d$. By definition 
$p(V(\sigma_F))= {\Delta}_{(d-1,1),F}$. As before, since the bundle 
${\mathcal F}_d$ is generated by the sections $\sigma_F$, we have that
$V(\sigma_F)$ is smooth of the right dimension, for a general $F$, 
and it is easy to verify that 
$p:V(\sigma_F) \to {\Delta}_{(d-1,1),F}$ is a desingularization.
Then we define $\tilde{\Delta}_{(d-1,1),F}:=V(\sigma_F)$.
Hence, using the adjunction formula, we can compute the 
canonical bundle of $\tilde{\Delta}_{(d-1,1),F}$ as 
the restriction to $\tilde{\Delta}_{(d-1,1),F}$ of the following
line bundle:
\begin{equation}\label{can2}
 K_{Gr(1,n)} + c_1({\mathcal F}_d)= 
 (2d-4)H + ( {{d(d+1)}\over {2}} - n-1 -2(d-1)) L.
\end{equation}

Consider  the bundles ${\mathcal N}^{d}_{\Gamma}$ and 
${\mathcal M}^{d}_{\Gamma}$ 
on $\Gamma$ respectively 
defined by the two following exact sequences :
\begin{eqnarray}\label{N_Gamma}
 0\rightarrow {\mathcal N}^{d}_{\Gamma}
 \rightarrow S^d \otimes {\mathcal O}_{\Gamma}
 \rightarrow {\mathcal F}_{d}\rightarrow 0,
\end{eqnarray}

\begin{eqnarray}\label{M_Gamma}
 0\rightarrow {\mathcal M}^{d}_{\Gamma}
 \rightarrow S^d \otimes {\mathcal O}_{\Gamma}
 \rightarrow {\mathcal E}_{d}\rightarrow 0.
\end{eqnarray}

From the definitions it follows that we have
\begin{eqnarray}\label{inclusion2}
 0\rightarrow {\mathcal M}^{d}_{\Gamma}
 \rightarrow {\mathcal N}^{d}_{\Gamma}
 \rightarrow {\mathcal K}\rightarrow 0.
\end{eqnarray}

The positivity result we will need this time is the following
\begin{lem}\label{pos3}
 If 
 \begin{equation}\label{cond3}
  {{d(d+1)}\over {2}} - n-1 -2(d-1)\geq c -1, 
 \end{equation}
 then  
 the bundle 
 ${\wedge^c{\mathcal M}^{d}_{\Gamma}}_{|{\tilde\Delta}_{(d-1,1),F}}\otimes 
 K_{{\tilde\Delta}_{(d-1,1),F}} (-H)$ is generated by its global sections.
\end{lem}
\begin{proof}
Use (\ref{can2}) and lemma \ref{pos}, (ii).
\end{proof}
Let now $\Delta_{d-1,1}\subset {\mathbf P}^n\times S^d$ be the family of the 
$\Delta_{(d-1,1),F}$'s,
and $\tilde\Delta_{d-1,1}\subset {\Gamma}\times S^d$ the family of the 
desingularizations. Let ${\cal Y}\subset \tilde\Delta_{r,d-r}$ be a subscheme
of relative dimension $k$, invariant under the action of
$GL(n+1)$, and $\tilde {\cal Y}\rightarrow {\cal Y}$ a desingularization.
Consider the sheaf $T^{vert}_{\tilde \Delta_{d-1,1}}$ defined by
$$
 0\rightarrow T^{vert}_{\tilde \Delta_{d-1,1}}
 \rightarrow T_{\tilde \Delta_{r,d-r}}
 {\rightarrow} T_{\Gamma}
 \rightarrow 0,
$$
and remark that
its restriction to ${{\tilde\Delta}_{(d-1,1),F}}$ coincides with 
${{\mathcal N}^{d}_{\Gamma}}_{|{\tilde \Delta_{(d-1,1),F}}}$.
If we assume $h^0 ({\tilde Y}_F,K_{{\tilde Y}_F} (-H))=0$ and set $c=n-1-k$,
we have that the natural adjunction map, 
\begin{equation}\label{lastvanishing}
 H^0 ({\wedge^c{\mathcal M}^{d}_{\Gamma}}_{|{\tilde\Delta}_{(d-1,1),F}}\otimes 
 K_{{\tilde\Delta}_{(d-1,1),F}} (-H))\to
 H^0({\tilde Y}_F,K_{{\tilde Y}_F}(-H)),
\end{equation}
that we can construct thanks to (\ref{inclusion2}), is obviously zero.

The last step will be the proof of the following
\begin{prop}\label{case2}
 Let $F$ be a general polynomial of degree $d$ verifying (\ref{cond3}).
 Let $Y_F\subset {\tilde\Delta}_{(d-1,1),F}$ 
 be a subvariety of dimension $k$, and $j:{\tilde Y}_F\to Y_F$ 
 a desingularization such that (\ref{lastvanishing}) vanishes.
 Then $Y_F$ is contained in the locus of lines of $X_F$.
\end{prop}
\begin{proof}
Recall that $F\in T^{vert}{\mathcal Y}_{|(x,\ell,F)}$. We claim 
that 
$$
 F\in {{\mathcal M}^{d}_{\Gamma}}_{|(x,\ell)}.
$$
Indeed if $F\in {\mathcal K}_{|(x,\ell)}$, then we have the surjection
\begin{equation}\label{ultsurj}
 T^{vert}{\mathcal Y}_{|(x,\ell,F)}\twoheadrightarrow 
 {\mathcal K}_{|(x,\ell)}.
\end{equation}
(This follows from the fact that if $F\in {\mathcal K}_{|(x,\ell)}$
then 
$$
 <S^1_x \cdot J_F^{d-1}, F>\twoheadrightarrow {\mathcal K}_{|(x,\ell)},
$$ 
plus Remark 2.0.1, (ii)).
Then, by (\ref{ultsurj}), we have that
$$
 codim_{{{\mathcal M}^{d}_{\Gamma}}{_|(x,\ell)}}
 T^{vert}{\mathcal Y}_{|(x,\ell,F)} = codim_{X_F} Y_F=c.
$$
As in Proposition \ref{case1}, we can now
use the commutative diagram
\begin{equation}\label{finaldiagram2}
  \xymatrix{
   H^0 (\wedge^c{{\mathcal M}^{d}_{\Gamma}}_{|{\tilde\Delta}_{F}}\otimes 
   K_{{\tilde\Delta}_{F}} (-H))
   \ar[d]^{ev\ \ }\ar@^{(->}[r]\ar@_{>}[r]
   &H^0 (\wedge^c{{T{\tilde \Delta}}}_{|{{\tilde\Delta}_{F}}}
   \otimes K_{\tilde\Delta_{F}}(-H))\ar[d]^{ev\ \ }\\
   (\wedge^c{{\mathcal M}^d_{\Gamma}}_{|{\tilde\Delta}_{F}}\otimes 
   K_{{\tilde\Delta}_{F}} (-H))|_{(x,\ell)}
   \ar[dr]^{<\cdot,{\bar W}>}\ar@^{(->}[r]\ar@_{>}[r]
   &(\wedge^c{{T{\tilde \Delta}}}_{|{{\tilde\Delta}_{F}}}
   \otimes K_{\tilde\Delta_{F}}(-H))|_{(x,\ell)}
   \ar[d]^{<\cdot,W>}\\
   &{\mathbb C}\\
}\end{equation}
and deduce from it, together with
Lemma \ref{pos3} and Remark 2.0.1, (ii), 
that what we claim holds, {\it i.e.}
$$
 F_{|\ell}=0.
$$
\end{proof}
The numerical condition (\ref{cond3}), which, for $c=n-1-k$, becomes
$$
  {{d(d-3)}\over {2}}\geq 2n-k-3,
$$
implies (\ref{cond2}). Thus, combining the two
propositions \ref{case1} and \ref{case2}, the proof of our main theorem is
completed.
\hfill $\Box$

\noindent Gianluca PACIENZA
\\Institut de Math\'ematiques de Jussieu
\\Universit\'e Pierre et Marie Curie
\\4, Place Jussieu, F-75252 Paris CEDEX 05 - FRANCE
\\e-mail: pacienza@math.jussieu.fr 

\smallskip

\noindent {\it current address:}
\\Department of Mathematics
\\Ohio State University
\\100 Mathematics Building
\\231 West 18th Avenue
\\Columbus, OH 43210-1174 - U.S.A.
\\e-mail: pacienza@math.ohio-state.edu

\end{document}